\documentstyle[txmac,a4,
amssymb,%
twoside,%
nocaphead,%
epsf,%
mypic,times,mathptm]{article}

\advance\oddsidemargin by -1.8cm
\advance\evensidemargin by -1.8cm
\advance\textwidth by 3.6cm

\def\mynewtheo#1#2{%
\newtheorem{@#1}{#2}
\newenvironment{#1}{\begin{@#1}\rm}{\end{@#1}}}

\mynewtheo{lemma}{Lemma}
\mynewtheo{theo}{Theorem}
\mynewtheo{rem}{Remark}
\mynewtheo{defi}{Definition}
\mynewtheo{conj}{Conjecture}
\mynewtheo{corollary}{Corollary}
\mynewtheo{proposition}{Proposition}
\mynewtheo{question}{Question}
\mynewtheo{exam}{Example}

\newenvironment{theorem}{\begin{theo}}{\end{theo}}

\newenvironment{eqn}{\begin{equation}}{\end{equation}}

\def\vrt#1{{\picfillgraycol{0}\picfilledcircle{#1}{0.09}{}}}
\def\fdot#1{{\picfillgraycol{0}\picfilledcircle{#1}{0.02}{}}}
\def\cycl#1#2#3#4{\vrt{#1}\vrt{#2}\vrt{#3}\vrt{#4}
\picline{#1}{#2}\picline{#2}{#3}\picline{#3}{#4}\picline{#4}{#1}}
\def\xcycl#1#2#3#4#5#6{\vrt{#1}\vrt{#2}\vrt{#3}\vrt{#4}\vrt{#5}
\vrt{#6}\picline{#1}{#2}\picline{#2}{#3}\picline{#3}{#4}
\picline{#4}{#5}\picline{#5}{#6}\picline{#6}{#1}}

\begin{document}

\makeatletter

\def\bysame{\same[\kern2cm]\,}
\def\lfra{\leftrightarrow}
\def\qed{\hfill\@mt{\Box}}
\def\@mt#1{\ifmmode#1\else$#1$\fi}
\def\qqed{\hfill\@mt{\Box\enspace\Box}}

\let\ap\alpha
\let\tl\tilde
\let\sg\sigma
\let\dl\delta
\let\Dl\Delta
\let\nb\nabla
\let\eps\varepsilon
\let\gm\gamma
\let\bt\beta
\def\cP{{\cal P}}
\def\bZ{{\Bbb Z}}
\def\bN{{\Bbb N}}
\def\len{\mbox{\operator@font len}\,}
\let\ds\displaystyle
\def\cf{\text{\rm cf}\,}
\def\md{\min\deg}
\def\Md{\max\deg}
\def\Mcf{\max\cf}
\def\br#1{\left\lfloor#1\right\rfloor}

\def\epsfs#1#2{{\catcode`\_=11\relax\ifautoepsf\unitxsize#1\relax\else
\epsfxsize#1\relax\fi\epsffile{#2.eps}}}
\def\epsfsv#1#2{{\vcbox{\epsfs{#1}{#2}}}}
\def\vcbox#1{\setbox\@tempboxa=\hbox{#1}\parbox{\wd\@tempboxa}{\box
  \@tempboxa}}
\def\p{\epsfsv{2cm}}

\def\@test#1#2#3#4{%
  \let\@tempa\go@
  \@tempdima#1\relax\@tempdimb#3\@tempdima\relax\@tempdima#4\unitxsize\relax
  \ifdim \@tempdimb>\z@\relax
    \ifdim \@tempdimb<#2%
      \def\@tempa{\@test{#1}{#2}}%
    \fi
  \fi
  \@tempa
}

\def\so{\Longrightarrow}

\def\go@#1\@end{}
\newdimen\unitxsize
\newif\ifautoepsf\autoepsftrue

\unitxsize4cm\relax
\def\epsfsize#1#2{\epsfxsize\relax\ifautoepsf
  {\@test{#1}{#2}{0.1 }{4   }
		{0.2 }{3   }
		{0.3 }{2   }
		{0.4 }{1.7 }
		{0.5 }{1.5 }
		{0.6 }{1.4 }
		{0.7 }{1.3 }
		{0.8 }{1.2 }
		{0.9 }{1.1 }
		{1.1 }{1.  }
		{1.2 }{0.9 }
		{1.4 }{0.8 }
		{1.6 }{0.75}
		{2.  }{0.7 }
		{2.25}{0.6 }
		{3   }{0.55}
		{5   }{0.5 }
		{10  }{0.33}
		{-1  }{0.25}\@end
		\ea}\ea\epsfxsize\the\@tempdima\relax
		\fi
		}

\let\old@tl\~
\def\~{\raisebox{-0.8ex}{\tt\old@tl{}}}

\author{A. Stoimenow\footnotemark[1]\ \,\footnotemark[2]\\[2mm]
\small Department of Mathematics, \\
\small University of Toronto,\\
\small Canada M5S 3G3\\
\small e-mail: {\tt alex@mpim-bonn.mpg.de},\\
\small WWW: {\hbox{\tt http://guests.mpim-bonn.mpg.de/alex}}
}

{\def\thefootnote{\fnsymbol{footnote}}
\footnotetext[1]{Supported by a DFG postdoc grant.}
\footnotetext[2]{On leave from: Max-Planck-Institut f\"ur Mathematik,
Vivatsgasse 7, D-53111 Bonn, Germany}
}

\title{\large\bf \uppercase{The skein polynomial of closed 3-braids}
\\[4mm]
{\small\it This is a preprint. I would be grateful
for any comments and corrections!}}

\date{\large Current version: \curv\ \ \ First version:
\makedate{28}{2}{2001}}

\maketitle

\long\def\@makecaption#1#2{%
   \vskip 10pt
   {\let\label\@gobble
   \let\ignorespaces\@empty
   \xdef\@tempt{#2}%
   }%
   \ea\@ifempty\ea{\@tempt}{%
   \setbox\@tempboxa\hbox{%
      \fignr#1#2}%
      }{%
   \setbox\@tempboxa\hbox{%
      {\fignr#1:}\capt\ #2}%
      }%
   \ifdim \wd\@tempboxa >\captionwidth {%
      \rightskip=\@captionmargin\leftskip=\@captionmargin
      \unhbox\@tempboxa\par}%
   \else
      \hbox to\captionwidth{\hfil\box\@tempboxa\hfil}%
   \fi}%
\def\fignr{\small\sffamily\bfseries}%
\def\capt{\small\sffamily}%

\newdimen\@captionmargin\@captionmargin2cm\relax
\newdimen\captionwidth\captionwidth\hsize\relax

\let\reference\ref
\def\eqref#1{(\protect\ref{#1})}

\def\proof{\@ifnextchar[{\@proof}{\@proof[\unskip]}}
\def\@proof[#1]{\noindent{\bf Proof #1.}\enspace}

\def\myfrac#1#2{\raisebox{0.2em}{\small$#1$}\!/\!\raisebox{-0.2em}{\small$#2$}}
\def\abstractname{}

\@addtoreset {footnote}{page}

\renewcommand{\section}{%
   \@startsection
         {section}{1}{\z@}{-1.5ex \@plus -1ex \@minus -.2ex}%
               {1ex \@plus.2ex}{\large\bf}%
}
\renewcommand{\@seccntformat}[1]{\csname the#1\endcsname .
\quad}

{\let\@noitemerr\relax
\vskip-2.7em\kern0pt\begin{abstract}
\noindent{\bf Abstract.}\enspace
Using the band representation of the 3-strand braid group,
it is shown that the genus of 3-braid links can be read off
their skein polynomial. Some applications are given, in particular
a simple proof of Morton's conjectured inequality and a condition
to decide that some polynomials, like the one of $9_{49}$, are
not admitted by 3-braid links. Finally, alternating links of
braid index 3 are classified.
\end{abstract}
}


\section{Introduction and results}

Braids are algebraic objects with a variety of applications.
They were defined and studied by Artin \cite{Artin,Artin2} and
although the connection to knot theory was known by Alexander
\cite{Alexander} and Markov \cite{Markov}, their importance in
this context was not recognized until the mid 80's. On the
topological side they were studied by Bennequin \cite{Bennequin}
in contact geometry and on the algebraic side used to discover
the Jones polynomial \cite{Jones} $V$ and its generalization, the
skein (HOMFLY) polynomial $P$ \cite{HOMFLY,PrzTr}, and some
relations between latter and braid representations were found,
as the Morton--Franks--Williams inequality \cite{Morton,WilFr}
(henceforth called MWF).

We will use henceforth the variable convention of \cite{Morton}
for $P$: $v^{-1}P(L_+)-vP(L_-)=zP(L_0)$. Here as usual
$L_{\pm,0}$ denote links with diagrams equal except near
one crossing, which is resp. positive, negative and smoothed out.

In \cite{Rudolph3}, Rudolph studied the topology of knots and links
via (braid closures of) representations of the braid
groups $B_n$ by embedded bands. These ``band'' representations
have been algebraically studied recently in \cite{BKL}, and before in
the 3-strand case by Xu.
%
%
In \cite{Xu}, he considered the new generator $\sg_3=\sg_1^{-1}
\sg_2\sg_1$ of $B_3$ ($\sg_{1,2}$ denoting Artin's generators),
with which $B_3$ gains the representation
\[
B_3\,:=\,\left\langle\,\sg_1,\,\sg_2\,\sg_3\,\big|\,
\sg_2\sg_1=\sg_3\sg_2=\sg_1\sg_3\,\right\rangle\,.
\]
As in \cite{Xu}, we set $\sg_{i\pm 3}=\sg_i$ (i.e., consider the
subscript only $\bmod 3$) to avoid awkward notation, in which case the
relations read $\sg_i\sg_{i-1}=\sg_{i+1}\sg_i$. We will also
sometimes denote a word $\sg_{s_1}^{\eps_1}\cdots \sg_{s_k}^{\eps_k}$
($\eps_i=\pm 1$) by $[l_1\dots l_k]$, where $l_i=\eps_is_i$.

Then for any representation of a braid $\bt\in B_3$ as word in
$\sg_{1,2,3}$ one obtains a Seifert surface of the closure
$\hat\bt$ of $\bt$ by inserting disks for each braid strand
and connecting them by half-twisted bands (the ``embedded bands'' in
Rudolph's terminology) along each $\sg_i^{\pm 1}$.
We will call this construction \em{band algorithm}; it has obviously
generalizations to higher braid groups. The surfaces thus obtained
were studied by Bennequin, who showed in particular:

\begin{theorem}(Bennequin \cite{Bennequin})\label{thB}
For every 3-braid link $L$ there exists a band (algorithm) Seifert
surface for $L$ of maximal Euler characteristic $\chi(L)$.
\end{theorem}

Based on this, Xu gave an algorithm to obtain the shortest
word representation of any $\bt\in B_3$ (in $\sg_{1,2,3}$) and thus
to calculate $\chi(\hat\bt)$. 

In this paper, we link his method with the skein polynomial and show

\begin{theorem}\label{th1}
For every 3-braid link $L$ we have $\Md_zP(L)=1-\chi(L)$, where
$\chi(L)$ is the maximal Euler characteristic of all Seifert
surfaces for $L$.
\end{theorem}

Thus we can relate for any 3-braid \em{knot} $K$ the two inequalities
(latter coming from \cite{Morton})
\[
\Md\Dl(K)\le g(K)=\frac{1}{2}\Md_zP(K)\le \tl g(K)\,,
\]
where $\Dl(K)$ is the Alexander polynomial \cite{Alexander2}, 
$g(K)$ is the genus of $K$ and $\tl g(K)$ its canonical genus, 
i.e. the minimal genus of the canonical Seifert surfaces of all
diagrams of $K$ (see e.g. \cite{gen1,Kobayashi}). There are
examples, like the knot $12_{2038}$ of \cite{KnotScape} (a picture
may be found in \cite{pos}), where the first (and very classical, see
e.g. \cite[exercise 10, p.\ 208]{Rolfsen}) inequality is not exact,
thus making the relation between the genus and $P$ even more
surprising (for this example we have $\Md\Dl=2$ and
$g=\myfrac{1}{2}\Md_zP=3$). Such examples do not occur among the
homogeneous knots of Cromwell \cite{Cromwell}.

Here are some straightforward and useful consequences of theorem
\reference{th1}:

\begin{corollary}\label{cr1}
There is no non-trivial 3-braid link with the skein polynomial
of the ($1,2,3$-component) unlink(s), and there are only finitely
many 3-braid links with the same skein polynomial.
\end{corollary}

This in particular allows to decide for a given $P$ polynomial
whether it is the polynomial of a 3-braid link. The general
problem whether the MWF bound can be realized among knots of given
$P$ polynomial was raised by Birman in problem 10.1 of \cite{Morton4}.
The negative answer was given in \cite{posex_bcr} by means of a
computer example and a braid index inequality of Jones \cite{Jones2},
but the problem for specific examples of polynomials, like those
of $9_{42}$ and $9_{49}$ (see problem 10.2 of \cite{Morton4}),
remained open. Now we can answer Birman's question negatively also
for these two polynomials. By Bennequin's theorem all 3-braid knots
of genus 2 can be easily written down (in fact all they have at
most 10 crossings, as will follow from an inequality proved below
in proposition \reference{pp1}, and hence are listed in the tables of
\cite{Rolfsen}), and no one has such a polynomial. (We will later
see that in fact for $9_{49}$ there is an even faster method to
conclude this, just looking at the polynomial.)

It also shows that, inspite of Birman's examples \cite{Birman},
later extended by Kanenobu \cite{Kanenobu}, the failure of
the skein polynomial to distinguish 3-braid links is limited,
and series of the type of \cite{Kanenobu2} do not exist among
such links.

\begin{corollary}\label{cr2}
For any 3-braid link $L$, we have $\md_vP(L)\le 1-\chi(L)$.
\end{corollary}

\proof The identity \cite[proposition 21]{LickMil} implies that
$\md_vP(L)\le \Md_zP(L)$ for any link $L$. (I'm grateful
to H.~Morton for pointing this out to me.)\qed

Corollary \reference{cr2} is a recent result of Dasbach--Mangum
\cite{DasMan}, a special case of a long-standing problem of
Morton \cite{Morton4}. Their proof is slightly less involved than
ours of theorem \reference{th1}, but the `$\le$' part of 
our equality follows rather easily (see lemma \reference{lm1})
so that we simplify the Dasbach--Mangum proof (in particular we
will not need the Scharlemann--Thompson result \cite{ST}).

\begin{rem}
In \cite{Morton}, Morton remarks that some knots
with $2g<\Md_zP$ exist. Since one of the 11 crossing knots with
unit Alexander polynomial, $11_{409}$ (denoted accoring to
\cite{KnotScape}), is such a knot (it has $g=2$ by \cite[figure
5 below]{Gabai}, $\Md_zP=6$ by \cite[p. 111]{LickMil}
and braid index 4, see \cite[figure 26]{LickMil}),
the `$\le$' inequality in theorem \reference{th1} is not true
for 4-braids. The problem of an example of a knot with the
contrary inequality $2g>\Md_zP$ is less straightforward to solve,
and an example was given only in \cite{posex_bcr}.
I have not investigated whether 4-braid examples exist.
\end{rem}

\section{The proof of theorem \reference{th1}}

We start with the easy part of theorem \reference{th1}, the
inequality `$\le $'.

\begin{lemma}\label{lm1}
For every 3-braid link $L$ we have $\Md_zP(L)\le 1-\chi(L)$.
\end{lemma}

\proof Because of theorem \reference{thB}, it suffices to prove
that for $L=\hat\bt$ we have
\begin{eqn}\label{!}
\Md_zP(L)\le \len(\bt)-2\,,
\end{eqn}
where $\len(\bt)$ is the length of $\bt$ as word in $\sg_{1,2,3}$.

We proceed by induction on $\len(\bt)$ (outer induction) and for
fixed value of $\len(\bt)$ on the crossing number $c(\bt)$ of $\bt$,
or equivalently by the number of letters $\sg_3^{\pm1}$ (inner
induction). In case $\sg_3^{\pm1}$ does not occur or $\len(\bt)\le 1$,
the inequality follows from \cite{Morton}, or is straightforward to
verify. Otherwise consider a letter $\sg_3^{\pm1}$ in $\bt$. If it
is followed by another letter $\sg_3^{\pm1}$, then we apply the skein
relation on one of the band crossings, use (outer) induction assumption
on $L_0$,
and switch the second band reverse to the first one, getting through
by outer induction. Else $\sg_3^{\pm1}$ is followed by $\sg_2^{\pm1}$
or $\sg_1^{\pm1}$. Again applying the skein relation on one of the
band crossings, we can switch it so that the subword of the 2 letters
reduces to one of the same length (two), but without occurrence of
$\sg_3^{\pm1}$, so we are done by the inner induction. \qed

\begin{rem}
Not only lemma \reference{lm1}, but in fact also the inequality
\eqref{!} (with `$2$' replaced by `$3$') is not true for 4-braids, since
the above quoted knot $11_{409}$ has a 7-band 4-braid representation:
$[(23-2)(1-2-1)2(2-3-2)(12-1)-32]$. (Here~-- and only here!~--
$\sg_3=[3]\in B_4$ does not refer to the element $[-121]\in B_3$.)
\end{rem}

The use of $\sg_{1,2,3}$ and Bennequin's result can be applied also for
the following useful inequality.

\begin{proposition}\label{pp1}
Let $c_3(K)$ be the minimal crossing number of a 3-braid representation
of a 3-braid knot or link $K$. Then
\begin{eqn}\label{q1}
c_3(K)\,\le\,\myfrac{5}{3}\bigl(3-\chi(K)\bigr)\,.
\end{eqn}
In particular if $K$ is a knot
\begin{eqn}\label{q2}
c_3(K)\,\le\,2\,\br{\frac{5}{3}\bigl(g(K)+1\bigr)}\qquad\mbox{or}
\qquad g(K)\ge \frac{3}{10}c_3(K)-1\,.
\end{eqn}
\end{proposition}

\proof Consider a minimal length (in $\sg_{1,2,3}$) word
representation for $\bt\in B_3$ with $\hat\bt=K$, and among such
word representations one of minimal crossing number (that is, minimal
number of letters $\sg_3^{\pm1}$), up to cyclic permutations.
Every (maximal) subword $\sg_3^{k}$, $k\in\bZ\setminus\{0\}$ of $\bt$ must
be (cyclically) followed by $\sg_2^{-1}$ or $\sg_1$ and preceded
by $\sg_2$ or $\sg_1^{-1}$, otherwise the first or last copy of
$\sg_3^{\pm1}$ can be eliminated by a relation, preserving the
word length of $\bt$. Since $\sg_3^{k}=\sg_1^{-1}\sg_2^{k}\sg_1$,
from each subword $\sg_3^{k}$ only one copy of $\sg_3^{\pm1}$
contributes three to the minimal crossing number of a
($\sg_{1,2}$-word) representation of $\bt$ (the others contribute $1$),
\eqref{q1} follows from Bennequin's result. For knots this is
equivalent to the second inequality of \eqref{q2}, the first inequality
follows by the remark that $2\mid c_3(K)$ for knots $K$. \qed

For $g=1$ we get from \eqref{q2} $c\le 6$, thus obtaining
Xu's list of $3_1$, $4_1$ and $5_2$ as the only genus 1 3-braid knots.
For $g=2$ we get $c\le 10$. The knots are in fact $3_1\#(!)3_1$,
$5_1$, $6_2$, $6_3$, $7_3$, $7_5$, $8_{20}$ and $8_{21}$
(compare the discussion after corollary \reference{cr1}).
There are three 12 crossing knots of genus 3, and still one
16 crossing knot of genus 4, so that \eqref{q2} is exact in
these cases.

We obtain from theorem \reference{th1} as corollary:

\begin{corollary}
If for the crossing number $c(K)$ of a knot $K$ it holds
\[
c(K)\,>\,2\,\br{\frac{5}{6}\bigl(\Md_zP(K)+2\bigr)}\,,
\]
then $K$ is not a closed $3$-braid.
\end{corollary}

It can already be expected from the proof of lemma \reference{lm1}
that the inequality $\Md_zP(L)\le 1-\chi(L)$ should be fairly
sharp. However, by computer check it turned out to be sharp without
any exception up to 18 bands, thus leading me to the investigation of
theorem \reference{th1}.

To carry out the rest of the proof of theorem \reference{th1},
we need to recall some of the work in \cite{Xu}. There a fast
algorithm to get any $\sg_{1,2,3}$ word-representation of $\bt\in
B_3$ into one of minimal length (and thus to calculate $\chi(\hat\bt)$)
is given.

We recall this algorithm as it will be important in the proof.

\def\labelenumi{\theenumi}
\def\theenumi{(\roman{enumi})}
\begin{enumerate}
\item Move all $\sg_i^{-1}$ to the left using $\sg_i\sg_j^{-1}=
\sg_{i+1}^{-1}\sg_{j+1}$. Thus $\bt=L^{-1}R$ with $L$ and $R$ positive.
\item As long as $L$ or $R$ contain some subword $\sg_{i+1}\sg_i$, this
subword can be moved to their beginning, giving $\bt=L^{-1}[21]^kR$
with $L$ and $R$ positive and non-decreasing.
\item Applying $[21]^{-1}\sg_{i+1}=\sg_i^{-1}$ and $\sg_i^{-1}[21]=
\sg_{i-1}$ and cyclic reductions, one of the 3 factors in
$L^{-1}[21]^kR$ can be eliminated.
\end{enumerate}

Although this may not be evident from the algorithm, we remark that
Bennquin's result implies that the minimal length of $\bt\in B_3$
is conjugacy invariant, and thus whether a word can be reduced by
Xu's algorithm or not is invariant under cyclic permutations of its
letters.

\proof[of theorem \reference{th1}] Since we need to consider only
one of two mirror images for $L$, we may assume at one point
in our proof for every case that $\bt\in B_3$
with $\hat\bt=L$ has non-negative exponent sum $e(\bt)$.

By Xu's algorithm, each $\bt\in B_3$ can be written in one of the
two forms
\def\theenumi{(\Alph{enumi})}
\begin{enumerate}
\item\label{itA} $[21]^kR$ or $L^{-1}[21]^{-k}$ ($k\ge 0$), or
\item\label{itB} $L^{-1}R$,
\end{enumerate}
where $L$ and $R$ are positive words with (cyclically)
non-decreasing indices (i.e. each $\sg_i$ is followed by $\sg_i$
or $\sg_{i+1}$). Since the form \reference{itB} must be cyclically
reduced, we may assume that $L$ and $R$ do not start or end with
the same letter.

If $\bt$ is of type \reference{itA} (we call this case ``strongly
quasipositive'' conforming to Rudolph \cite{Rudolph,Rudolph2}), then by
the mirroring argument we may assume $e(\bt)>0$, and then have
from \cite{Morton} and theorem \reference{thB}
\begin{eqn}\label{q}
1-\chi(\hat\bt)\,=\,e(\bt)-2\,\le\,\md_vP\,\le\,\Md_zP\,,
\end{eqn}
and thus the reverse inequality to lemma \reference{lm1}. 

Thus we need to consider only the case \reference{itB}.

A fair part of our argument will go like this: We choose a band
(crossing) in $\bt=\bt_{\pm}$ and apply the skein relation at this
crossing, expressing the polynomial of $L=L_{\pm}=\hat\bt$ by those of
$L_0=\hat\bt_0$ and $L_{\mp}=\hat\bt_{\mp}$. (Here $\bt_0$ and
$\bt_{\mp}$ are obtained by deleting resp. reversing the
band in $\bt$ we consider.) Then we show that only one of $\bt_0$
and $\bt_{\mp}$ contributes to the coefficient $[P(L)]_{z^{1-\chi(L)}}$
of $z^{1-\chi(L)}$ in $P(L)$. Because of lemma \reference{lm1} for
this it suffices to show that the other one is not of Xu's minimal
(word length) types \reference{itA} and \reference{itB}.
This way we lead back inductively the case of $L$ to some simple cases.

We have
\begin{eqn}\label{q3}
\bt=L^{-1}R \mbox{\quad with\quad } R=\prod_{j=0}^l\sg_{i+j-1}^{k_j}
\mbox{\quad and\quad } L=\prod_{j=0}^{l'}\sg_{i'+j-1}^{k'_j}\,,
\end{eqn}
with $i\ne i'$ and $i+l\not\equiv i'+l'\bmod 3$.

The first application of the skein relation argument is that
we can make induction on the $k_j$ and $k'_j$, thus being
left just with the cases where all $k_j=k'_j=1$, in which case $L$ and
$R$ get the simpler form
\begin{eqn}\label{q4}
R=(\sg_{i}\sg_{i+1}\sg_{i+2})^k\ap\mbox{\quad and\quad }
L=(\sg_{i'}\sg_{i'+1}\sg_{i'+2})^{k'}\ap'
\end{eqn}
with $\ap$ and $\ap'$ of length $\le 2$. Again we can assume
modulo mirror images, that $R$ is not shorter than $L$, i.e. $\bt=
L^{-1}R$ with $L$ and $R$ as in \eqref{q4} has $e(\bt)\ge 0$ (it
may originate from a braid in \eqref{q3} with negative exponent sum!).

Now consider the case where $6\mid e(\bt)$ and use the representation
theory of $P$ on 3-braids (see \cite{Jones2}). Let $\Dl=[121]$ be the
square root of the
center generator of $B_3$ \cite{Chow}. Define $\bt^*=\Dl^{2/3e(\bt)}
\bt^{-1}$ to be the \em{dual} of $\bt$ (clearly $\bt^{**}=\bt$).
Then, as
observed in \cite[proof of proposition 2]{Birman}, $\hat\bt$ and
$\hat\bt^*$ have the same polynomial, because they have the same
(normalized) Burau trace and the same exponent sum. But for $\bt$
as in \eqref{q3} we have because of $\Dl^2\sg_{i+1}^{-1}\sg_i^{-1}
\sg_{i-1}^{-1}=\sg_i^3$ that
\begin{eqn}\label{q45}
\bt^*\,=\,\gm\sg_i^{k_i}\sg_j^{k_j}\gm'\,,
\end{eqn}
where $\gm$ and $\gm'$ have length at most 2, $i\ne j$ and
$k_ik_j\le 0$, and are thus left with showing that for such
words $\bt^*$, $\Md_zP(\hat\bt^*)=\len(\bt^*)-2$. Again by the
skein induction argument on the $k_i$ and $k_j$ this can be reduced
to the cases where $\bt^*$, and hence $\bt$, have small crossing
number, and they checked directly.

Now consider the case where $e(\bt)\equiv 1,2\bmod 6$. We apply the
skein relation at the rightmost letter/band in $\bt$. Then one of
$\bt_0$ or $\bt_-$ have $6\mid e$. It suffices to show that the other
one is not minimal (and apply lemma \reference{lm1}). For this one
checks that either $\bt_0$ is not cyclically reduced (starts and
ends with opposite letters), or that when permuting the
rightmost (negative) letter of $\bt_-$ to the left, the word $L$
in \eqref{q3} is not increasing (and hence $\bt_-$ can be reduced by
Xu's algorithm). For example for $\bt=[-2-1-3-21231231]$ we get
$\bt_-\doteq [-1-2-1-3-2123123]$ (`$\doteq$' meaning equality up to
conjugacy).

If $e(\bt)\equiv 3,4,5\bmod 6$, then apply the same argument at most
3 times, getting back to the $6\mid e$ case (except in the cases
where $R$, and hence $L$, are short, and which can be checked
directly). Since any $6\mid e$ word is reduced, and
every pair $(\bt_0,\bt_-)$ contains one reducible
word, it will indeed be the $6\mid e$ braid wo which the argument
recurs rather than some of its neighbors.

This completes the proof of theorem \reference{th1}. \qed

\begin{rem}
There is an alternative way to proceed with the
proof after \eqref{q4}, namely to remark that in the application of
the skein relation at every second stage it is $\bt_0$ that is of
Xu's form, and then to work by induction on the word length.
Thus the representation theoretic argument can be avoided.
However, the proof did not appear (to me) more elegant without it, and
also, there are some insights which this argument explains better
(in particular the cases of trivial Alexander polynomial, see
question \reference{qu1}), so I consider it not inappropriate.
\end{rem}

The representation theory also shows that $P(\hat\bt)$ for
$\bt\in B_3$ can be calculated in time $O(c(\bt))$ (see \cite{%
MorSho1}), and thus we have an even faster algorithm than the
(quadratic) one of Xu to calculate $\chi(\hat\bt)$.

\begin{corollary}
For $\bt\in B_3$, $\chi(\hat\bt)$ can be calculated in $O(c(\bt))$
steps.
\end{corollary}

\section{Further applications}

We can even say a little more that theorem \reference{th1}.
Since throughout the proof, $L_+$ inherited its maximal coefficient
from $L_-$ or $L_0$ (up to multiplication with units in
$\bZ[v,v^{-1}]$) and duality does not alter the polynomial, we
see that in fact we can determine what maximal ($z$-)coefficients
skein polynomials of $3$-braid links can have by checking some
simple cases. We have the following result (note that it
somewhat depends on the convention for $P$ chosen!):

\begin{theorem}\label{th2}
Let $L$ be a $3$-braid link. Then $[P(L)]_{z^{1-\chi(L)}}$ is up to
units $\pm v^k$, $k\in\bZ$ one of $1\pm v^2$ or $1$, except for the
3 component unlink (where it is $(1-v^2)^2$). If $L$ is a knot
or 3 component link, then $[P(L)]_{z^{1-\chi(L)}}\ne -v^k(1+v^2)$.
\end{theorem}

\proof[(sketch)] This is, as remarked, basically a repetition of the
proof of theorem \reference{th1}. In the strongly
quasipositive case (which there could be dealt with immediately)
the skein and duality arguments can be applied similarly, leaving us
with a braid of the form
$(\sg_2\sg_1)^k\sg_1^l\ap'$ with $\ap'$ having small length. By skein
argument induction, $l$ can be reduced to $1$. For the first factor
use now $(\sg_2\sg_1)^3=\sg_2\sg_1\sg_2^2\sg_1\sg_2$, and apply
the skein argument on the `$\sg_2^2$' in the middle until you
get $k$ small. The rest for the first statement is to compute
the polynomial for some simple words.

To show that $(1-v^2)^2$ occurs only for the 3-component unlink,
we need to verify this among the small words and to observe that
the procedure of inductively simplifying the braid in the
proof of theorem \reference{th1} at no stage gives the
trivial (empty) word (for this $i\ne j$ in \eqref{q45} is needed).

For the second statement, consider $\hat\bt=L$ for $\bt$ of odd
connectivity (i.e. even exponent sum) and the signs in the skein
relations.
It follows from the skein relation that expressing $P(L_\eps)$ for
$|\eps|=1$ by $P(L_{-\eps})$ and $P(L_0)$, latter's coefficients are
$+1$ except for the one of $P(L_0)$ when $\eps=-$. Thus we
need to take care only when we switch negative crossings.
When reducing the $k_j'$ in \eqref{q3}, we can maintain
sign at the cost of leaving possibly one of them equal to 2.
Denote by $w$ the subword of $L$ made up of this generator square.
Then, a possible mirroring (to get in \eqref{q4} $R$ to be not
shorter than $L$) does not alter the sign of $[P(L)]_{z^{1-\chi(L)}}$.
When mirroring puts $w$ into $R$, the generator square
has positive sign and can be reduced. Then, after going over to
$\bt^*$, (except for the few small length cases where $R$ is short) we
switch negative crossings only when reducing the negative one of $k_i$
and $k_j$ in \eqref{q45}. Then we just choose to reduce it by steps of
$2$, thus preserving connectivity of the closure and sign of
$[P(L)]_{z^{1-\chi(L)}}$. If $w$ remains in $L$ (in which case we don't
apply mirroring), reduce it as well, but then in $\bt^*$ in
\eqref{q45} reduce the negative one of $k_i$ and $k_j$ first by one,
thus canceling the negation. Checking some simple cases (just of odd
connectivity)
shows the result up to mirroring. Since mirroring does not negate 
$[P(L)]_{z^{1-\chi(L)}}$ when the connectivity is odd (or equivalently
$2\mid e(\bt)$), the result follows. \qed

We can say something on the cases where $1\pm v^2$ in the above
theorem occurs as maximal coefficient. We rephrase this using the
relation to the Conway polynomial $\nb(K)$ \cite{Conway} and
Alexander polynomial $\Dl(K)$ (in the normalization $\Dl(1)=1$ and
$\Dl(t^{-1})=\Dl(t)$)
\begin{eqn}\label{q6}
\Dl(t)\,=\,\nb(t^{1/2}-t^{-1/2})\,=\,P(1,t^{1/2}-t^{-1/2})\,.
\end{eqn}

\begin{proposition}\label{pMcf}
If for a 3-braid link $L$, $\Md\nb(L)<\Md_z P(L)$ or $\Mcf\nb(K)=\pm
2$, then $L$ is (the closure of) a strongly quasipositive 3-braid.
\end{proposition}

\proof Again check the small length cases and apply the previous
type of induction. \qed

A small application of this is

\begin{corollary}
Any homogeneous braid index 3 link $L$ is fibered or positive.
\end{corollary}

\proof Theorem \reference{th2} shows from \eqref{q6} that $\Dl(L)$ has
leading coefficient $\Mcf\Dl=-1$, $1$ or $2$, since $2\Md\Dl(L)=
\Md_zP(L)$ by \cite{Cromwell}. If the leading coefficient is $\pm 1$,
then $L$ is fibered (see \cite[corollary 5.3]{Cromwell}). Otherwise,
it is strongly quasipositive, and lemma \reference{lm1} and \eqref{q}
imply $\md_vP(L)=\Md_zP(L)$. Then apply \cite[theorem 4]{Cromwell}. \qed

We will in the next section have to say much more about
alternating links.

Some other worth remarking consequences follow now from the
work of Rudolph \cite{Rudolph2}. For simplicity, call the
Alexander polynomial $\Dl(K)$ of a knot $K$ \em{maximally
monic}, if its leading coefficient is $\pm1$ (monicness) and
its degree equal to $g(K)$ (maximality). A classical result states
that fibered knots have such Alexander polynomials. Here we obtain:

\begin{corollary}
Any achiral or slice braid index 3 knot has maximally monic
Alexander polynomial.
\end{corollary}

\proof For slice knots this follows from proposition \reference{pMcf}
and \cite{Rudolph2}. For achiral knots use theorem \reference{th2}
and that $[P(K)]_{z^{1-\chi(K)}}$ is self-conjugate. \qed

In the slice case neither maximality nor monicness need to hold for
4-braids, as show (slice) knots like $8_8$ and the $\Dl=1$ 11 crossing
knot $11_{409}$. In the achiral case the situation
is unclear since there may exist no achiral braid index 4 knot
(see \cite{posex_bcr}).

In a similar way we get

\begin{corollary}
There are only finitely many braid index 3 knots $K$ of given unknotting
number $u(K)$, whose Alexander polynomial is not maximally monic. (For
unknotting number 1 this is just the knot $5_2$.)
\end{corollary}

\proof Use that by \cite{Rudolph2}, $u(K)\ge g(K)$ for 
strongly quasipositive $K$. \qed

This is certainly not true, already for $u(K)=1$, without the condition
on the Alexander polynomial, an example being the rational knots with
Conway notation $(n11n)$, $n\in\bN$.

Finally we remark that theorem \reference{th2} also gives another (and
much more straightforward) way to see that $9_{49}$ (and any other knot
with such polynomial) is not a 3-braid knot. Unfortunately, this
simple criterion does not always work, as shows the polynomial of
$9_{42}$. This can always be decided as discussed after
corollary \reference{cr1}, but for higher $\Md_zP$ the process
of generating the whole list of knots becomes tedious, so that
our work here does not render obsolete examples like the one
in \cite{posex_bcr}.

\section{Alternating links of braid index 3}

A final, and main, application of our method is
to complete the description of alternating links of braid index 3.
Murasugi \cite{Murasugi} described the rational ones among them.
Our result easily implies his.

\begin{theorem}\label{tha3}
Ler $L$ be an alternating braid index 3 link.
Then (and only then) $L$ is
\def\theenumi{\alph{enumi})}
\begin{enumerate}
\item\label{scA} the connected sum of two $(2,k)$-torus links
(with parallel orientation), or
\item\label{scB} an alternating 3-braid link (i.e. the closure of an
alternating 3-braid, including split unions
of a $(2,k)$-torus link and an unknot and the 3 component unlink), or
\item\label{scC} a pretzel link $\cP(1,p,q,r)$ with $p,q,r\ge 1$
(oriented so that the twists corresponding to $p,q,r$ are
parallel).
\end{enumerate}
\end{theorem}

\proof
We know from \cite{Cromwell} that for an alternating link $L$,
$\Md\nb(L)=\Md_zP(L)$, and thus for braid index $b(L)\le 3$ we have
from proposition \reference{pMcf} that $\Mcf\nb(L)\in\{-1,1,2\}$.

In case $\Mcf\nb(L)=\pm 1$, $L$ is fibered by \cite{Murasugi2} (or see
\cite[corollary 5.3]{Cromwell}), and then by \cite[theorem A(2)]{Murasugi}
any alternating diagram of $L$ has $b(L)\le 3$ Seifert circles.
This gives the cases \reference{scA} and \reference{scB}. (The
split cases are easy since $b$ is additive under split uniton.) Since
it is known from \cite{Murasugi3} that case \reference{scA} includes
all composite links of braid index 3, we may henceforth assume
that $L$ be prime, and also non-split.

Assume now that $\Mcf\nb(L)=2$. We know from proposition
\reference{pMcf} that $L=\hat\bt$ with $\bt\in B_3$ strongly
quasipositive. By lemma \reference{lm1} and \eqref{q} we have $\md_vP
(L)=\Md_zP(L)$, and then it follows from \cite[theorem 4]{Cromwell}
(see also \cite{Nakamura}) that any homogeneous (in particular,
alternating) diagram of $L$ is positive. Thus $L$ has a special
alternating diagram $D$.

For every such diagram $D$ we consider the Seifert graph $G(D)$,
with vertices corrsponding to Seifert circles and egdes to
crossings (see \cite[\S 1]{Cromwell}). $G$ is connected, planar and
bipartite, hence every cycle in $G$ has even length (possibly
$2$, since $G$ may have multiple edges). For every such $G$
we can contrarily construct a special alternating diagram $D(G)$
with $G(D(G))=G$ (which depends on the planar embedding of $G$
only modulo flypes). It follows from \cite{Cromwell} (see corollary
2.2 and the proof of theorem 5, p.\ 543)
that if $G'$ is obtained from $G$ by deleting an edge, then
$\Mcf\nb(D(G'))\le \Mcf\nb(D(G))$. (Here deleting an edge $e$ means
deleting one single edge in a multiple one. If $e$ is single
and its deletion disconnects the graph, then if one of the 2
new components is a single vertex, this vertex is deleted as well,
while if both components contain edges, the deletion of $e$
is prohibited.)

If $G$ is a cycle graph of length $2k$ like
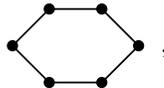
\[
\diag{7mm}{3}{1.4}{\xcycl{1 0}{2 0}{2.7 0.7}{2 1.4}{1 1.4}{0.3 0.7}}\,,
\]
then $D(G)$ depicts the $(2,2k)$-torus link $T_k$ with reverse
orientation, and $\nb(T_k)=kz$. Therefore, if $\Mcf\nb(D(G))=2$,
$G$ cannot contain a cycle of length $>4$.

Since we excluded composite links, we may assume that in $G=G(D)$
there is no vertex connected (by a possibly multiple edge)
to only one single other vertex. Then we replace in $G$
every multiple edge by a single one. We obtain a graph $\hat G$
(called sometimes reduced Seifert graph), in which each vertex has
valency $\ge 2$ and there are no multiple edges. By
\cite[proposition 13.25]{BurZie} (see also \cite{Cromwell,restr}),
$\Mcf\nb(D(\hat G))=1$ iff $\hat G$ is a tree, and in this case
$\hat G$ would have to be one single vertex, which is
uninteresting.

Therefore, $\Mcf\nb(D(\hat G))=2$ and $\hat G$ still contains a cycle.
We know from $G$ that any cycle in $\hat G$ has length $4$. We wish
to show that there is only one such cycle.

Assume there were two, call them $C_1$ and $C_2$. If $C_1$ and
$C_2$ have $\le 1$ vertex in common, then by deleting edges from
$\hat G$ we can obtain a graph $\tl G$ consisting of $C_1$ and $C_2$
joined by a (possibly trivial) path.
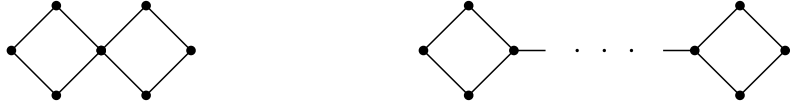
\[
\diag{6mm}{4}{2}{\picmultigraphics{2}{2 0}{\cycl{0 1}{1 0}{2 1}{1 2}}}
\kern3cm
\diag{6mm}{8}{2}{
  \picmultigraphics{2}{6 0}{\cycl{0 1}{1 0}{2 1}{1 2}}
  \picline{2 1}{2.7 1}\picline{6 1}{5.3 1}
  \picmultigraphics{3}{0.6 0}{\fdot{3.4 1}}
}
\]
But $D(\tl G)=T_2\# T_2$, and $\Mcf\nb(T_2\# T_2)=4$. If $C_1$ and $C_2$
have two neighbored vertices in common, then $G$ has a subgraph
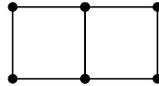
\[
\diag{6mm}{3.2}{2.0}{
  \pictranslate{0 0.4}{
    \picmultigraphics{2}{1.6 0}{\cycl{0 0}{0 1.6}{1.6 1.6}{1.6 0}}}
  }
\]
and a cycle of length $6$. Thus either $C_1$ and $C_2$ have 3 vertices
in common or two vertices which are opposite (not neighbored).
In both cases $G$ contains the subgraph
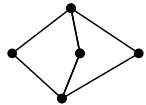
\[
\diag{6mm}{3}{2}{
  \cycl{1.3 2}{0 1}{1.1 0}{1.5 1}
  \cycl{1.3 2}{2.8 1}{1.1 0}{1.5 1}
}
\]
This corresponds to the $(2,2,2)$-pretzel link (oriented so that
the clasps are reverse) with $\nb=3z^2$.

Therefore, $\hat G$ contains only one cycle (of length 4), and must
be only this cycle. This shows that $D$ is a diagram of the
$(p,q,r,s)$-pretzel link (with parallel twists) $\cP(p,q,r,s)$.
Since for $p=1$ we have $\cP(1,q,r,s)=\widehat{\ }\,[1^q2^s3^r]$, it
remains to show that $L=\cP(p,q,r,s)$ has braid index $4$ for
$p,q,r,s>1$. Using MWF and $\md_vP(L)=\Md_zP(L)=1-\chi(L)$, it
suffices to show that $\mu(L):=\Md_vP(L)=7-\chi(L)$. For this, we
verify it for $p=q=r=s=2$ and inductively use the skein relation,
noticing that the signs of the $z$-coefficients of $[P(L')]_
{v^{\mu(L')}}$ are for $L'=L_-$ and $L'=L_0$ the same as for
$L'=\cP(2,2,2,2)$, and thus their contributions to $P(L_+)$ do not
cancel.  \qed

\begin{rem}
K.\ Murasugi pointed out that an alternative proof of the
conclusion $L$ special alternating and $\Mcf\Dl_L=2\ \so
L=\cP(p,q,r,s)$ was given in lemma 4.3 of \cite{Murasugi5}.
\end{rem}

By \cite{Kauffman,Murasugi,Thistle}, each alternating 3-braid knot
will have even crossing number. The theorem now shows:

\begin{corollary}
Prime alternating braid index 3 knots, which are
not closures of alternating 3-braids, have odd crossing number.
\end{corollary}

Also we have

\begin{corollary}
Each alternating braid index 3 link is an alternating 3-braid link
or is positive.
\end{corollary}

\begin{rem}
The braid representations of 3-braid links were described in
\cite{BirMen}, but since braids have (at least so far) proved
of little use in the study of combinatorial (diagrammatic) properties
of their closures, the methods there are unlikely to approach such
kind of results.
\end{rem}

\section{Problems}

Here are some open questions one can ask. For example, one is
the following question, 
suggested by computer experiment, in which braid index 3 knots
of at most 16 crossings were identified in the tables of
\cite{KnotScape} and all were found to accord to the following
conjectured rule (the same experiment pointed me to theorem
\reference{tha3}).

\begin{question}
Does any non-alternating braid index 3 knot have even crossing number?
\end{question}

Another problem which is possible to pursue by the methods of this
paper, but which involves some technical difficulties is

\begin{question}\label{qu1}
Does for any 3-braid link $K$ with $\nb(K)\ne 0$ hold $\Md\nb(K)\ge
\Md_z P(K)-2$? Are the only 3-braid links $K$ with $\nb(K)=0$
the split unions of $(2,k)$-torus links and an unknot and
links of the form $\hat{\ }\bigl([123]^{2k}\bigr)$?
\end{question}

The proof should go similarly to theorem \reference{th1}, but
more care must be taken.

The origin of the investigations of this paper came
from the attempt to compare the two estimates for $\tl g(K)$
given by $g(K)$ and $\myfrac{1}{2}\Md_zP(K)$. Now it was shown that
they are equally good, but I do not know whether they are always sharp.

\begin{question}
Is for any 3-braid knot $K$, $g(K)=\tl g(K)$?
\end{question}

Since, as mentioned, $2\tl g(K)=\Md_z P(K)$  for $K$ of $\le 12$
crossings, the answer is positive at least up to genus 3.

A final question concerns a possible generalization of
Bennequin's result.

\begin{question}\label{qu4}
Does any knot have a minimal genus Seifert surface constructed
by the band algorithm on a minimal strand representation?
\end{question}

It has been mentioned by Birman and Menasco that 
knots lacking the requested property 
should exist. However, I was unable to find (by computer) a concrete
example. In \cite{Rudolph3}, Rudolph showed that any Seifert surface is
isotopic to some band algorithm surface (for a possibly non-%
minimal strand representation).

This question will be answered in a joint paper with M.~Hirasawa
\cite{HS}.


\noindent{\bf Acknowledgement.} I would wish to thank to L.~Rudolph,
H.~Morton and K.~Murasugi for helpful remarks and discussions and
to M.~Hirasawa for his collboration in examining question
\reference{qu4}.

{\small
\let\old@bibitem\bibitem
\def\bibitem[#1]{\old@bibitem}

}
\end{document}